# COSET DIAGRAM FOR THE ACTION OF PICARD GROUP ON $\mathbb{Q}(i, \sqrt{3})$

**Qaiser Mushtaq** [a], **Saima Anis** [b,1]

[a] Department of Mathematics, Quaid-i-Azam University, Islamabad, Pakistan

qmushtaq@isb.apollo.net.pk

[b] Department of Mathematics, COMSATS Institute of Information Technology, University Road, Tobe Camp, Abbottabad, Pakistan

saimaanis_pk@yahoo.com

**Abstract.** The Picard group $\Gamma$ is $PSL(2, \mathbb{Z}[i])$. We have defined coset diagram for the Picard group. It has been observed that some elements of $\mathbb{Q}(i, \sqrt{3})$ of the form $\frac{a+b\sqrt{3}}{c}$ and their conjugates $\frac{a-b\sqrt{3}}{c}$ over $\mathbb{Q}(i)$ have different signs in the coset diagram for the action of $\Gamma$ on the biquadratic field $\mathbb{Q}(i, \sqrt{3})$, these are called ambiguous numbers. We have noticed that ambiguous numbers in the coset diagram for the action of $\Gamma$ on $\mathbb{Q}(i, \sqrt{3})$ form a unique pattern. It has been shown that there are finite number of ambiguous numbers in an orbit $\Gamma \alpha$, where $\alpha$ is ambiguous, and they form a closed path and it is the only closed path in the orbit $\Gamma \alpha$. We have devised a procedure to obtain ambiguous numbers of the form $\frac{a+k\sqrt{3}}{c}$, where $k$ is a positive integer.

**2000 Mathematics Subject Classification:** Primary 20G15, Secondary 11E04 and 05C25.

**Keywords:** Picard group, biquadratic field, coset diagram and ambiguous number.

## Introduction

Picard group $\Gamma$ or $PSL(2, \mathbb{Z}[i])$ is the group of linear fractional transformations $T(z) = \frac{az+b}{cz+d}$ with $ad - bc = 1$ and $a, b, c, d \in \mathbb{Z}[i]$.

The finite presentation of $\Gamma$ is
$\langle A, B, C, D : A^3 = B^2 = C^3 = D^2 = (AC)^2 = (AD)^2 = (BC)^2 = (BD)^2 = 1 \rangle$, where $A, B, C$ and $D$ are linear fractional transformations defined by $A(z) = \frac{1}{z-i}$, $B(z) = \frac{1}{z}$, $C(z) = \frac{1+z}{-z}$ and $D(z) = \frac{-1}{z}$, for details see [3], [4] and [10].

In 1988 [6] Q. Mushtaq has defined ambiguous numbers in real quadratic fields and studied these numbers for an action of the modular group on real quadratic fields $\mathbb{Q}(\sqrt{n})$. Later in 1995 [2] M. Aslam, S.M. Husnine and A. Majeed have found a formula which determines the exact number of ambiguous numbers in the coset diagrams for the action of the modular group on $\mathbb{Q}(\sqrt{n})$. Since Picard group is an extension of the modular group, therefore we have investigated ambiguous numbers in biquadratic field by the action of $\Gamma$ on the biquadratic field.

---

[1] This paper was written in 2007 after I have given seminar in Algebra Forum (www.algebraforum.org.pk/about/). It is from my Ph.D. thesis completed in July 2009. I have also given talk of some of its results in 10th International Pure Mathematics Conference in August 2009 in Islamabad, Pakistan.

The Picard group $\Gamma$ acts on $\mathbb{Q}(i)$ and its extensions. The fixed points of generators $A, B, C$ and $D$ of $\Gamma$ are $\frac{i \pm \sqrt{3}}{2}, \pm 1, \frac{-1 \pm \sqrt{3}i}{2}$ and $\pm i$ respectively. They all lie in a biquadratic field $\mathbb{Q}(i, \sqrt{3})$, where $i$ and $\sqrt{3}$ are zeros of an irreducible polynomial $(t^2 - 3)(t^2 + 1)$ over $\mathbb{Q}$. The action of $\Gamma$ on $\mathbb{Q}(i, \sqrt{3})$ is different from $\mathbb{Q}(i, \sqrt{n})$, where $n > 1$ and $n \neq 3$ is a square-free integer, and deserves special account. Therefore the closed paths in the coset diagram for the action of $\Gamma$ on $\mathbb{Q}(i, \sqrt{3})$ must be significantly different from the coset diagrams for the action of $\Gamma$ on $\mathbb{Q}(i, \sqrt{n})$ when $n \neq 3$.

The elements of $\mathbb{Q}(i, \sqrt{3})$ are of the form $u + v\sqrt{3}$, where $u, v \in \mathbb{Q}(i)$. They can be written as $\alpha = \frac{(a+bi)+(c+di)\sqrt{3}}{e}$, where $a, b, c, d, e \in \mathbb{Z}$. The conjugates of $\alpha$ over $\mathbb{Q}$ are $\alpha_1 = \frac{(a+bi)-(c+di)\sqrt{3}}{e}, \alpha_2 = \frac{(a-bi)+(c-di)\sqrt{3}}{e}$ and $\alpha_3 = \frac{(a-bi)-(c-di)\sqrt{3}}{e}$. The conjugate of $\alpha$ over $\mathbb{Q}(i)$ is $\alpha_1$ and the conjugate of $\alpha$ over $\mathbb{Q}(\sqrt{3})$ is $\alpha_2$. The action of $\Gamma$ on $\mathbb{Q}(i, \sqrt{3})$ shows that certain elements of $\mathbb{Q}(i, \sqrt{3})$ of the form $\frac{a+b\sqrt{3}}{c}$ behave special under this action. Thus they deserve a classification. As there are two conjugates of $\alpha = \frac{a+b\sqrt{3}}{c}$ over $\mathbb{Q}$, namely, $\alpha$ and $\frac{a-b\sqrt{3}}{c}$, and the conjugate of $\alpha$ is again $\alpha$ over $\mathbb{Q}(\sqrt{3})$, so we have considered conjugate of $\alpha$ over $\mathbb{Q}(i)$, that is, $\bar{\alpha} = \frac{a-b\sqrt{3}}{c}$. A real quadratic irrational number $\alpha = \frac{a+b\sqrt{3}}{c} \in \mathbb{Q}(i, \sqrt{3})$, where $a, b, c \in \mathbb{Z}$, is called totally positive (negative) if $\alpha$ and $\bar{\alpha}$ are both positive (negative). When $\alpha$ and $\bar{\alpha}$ have opposite signs, then they are called ambiguous numbers. A general element $\frac{(a+bi)+(c+di)\sqrt{3}}{e}$ of $\mathbb{Q}(i, \sqrt{3})$ does not have any sign that is why we have considered elements of the form $\frac{a+b\sqrt{3}}{c}$. We are discussing orbits of $\mathbb{Q}(i, \sqrt{3})$ by action of $\Gamma$ therefore we have written in the statements of our results that ambiguous numbers are in $\mathbb{Q}(i, \sqrt{3})$. They play an important role in classifying the orbits of $\mathbb{Q}(i, \sqrt{3})$ when $\Gamma$ acts on it. A diagrammatic argument, called coset diagram for the action of Picard group $\Gamma$ on $\mathbb{Q}(i, \sqrt{3})$, is defined and used to prove results in this paper. G. Higman and Q. Mushtaq have defined coset diagram for modular group in [5]. One of the applications of coset diagram for modular group is in [1]. We have defined coset diagram for the Picard group $\Gamma$. They need symbols for the generators as well as a method or pattern to join them. The group $\Gamma$ consists of four generators, two of order 3 and two of order 2, so it is possible to avoid using colours. The generators $A$ and $C$ both have order $3$, so the $3$-cycles of $A$ and $C$ are represented by triangles. But to distinguish generator $A$ from generator $C$, we have denoted the $3$-cycles of the generator $C$ by three unbroken edges of a triangle permuted anticlockwise. The $3$-cycles of the generator $A$ are denoted by three broken edges of a triangle permuted anticlockwise.

As generators $B$ and $D$ are involutions so we have represented them by edges and orientation of edges can be avoided. To distinguish generator $B$ from generator $D$, the 2-cycles of generator $B$ have been represented by a bold edge and two vertices which are interchanged by $D$ have been joined by an edge. Fixed points of $A, B, C$ and $D$, if they exist, have been denoted by heavy dots.

The fragment of a coset diagram in Fig.1 explains clearly the amalgam structure of $\Gamma$, that is, $\Gamma = (A_4 \underset{Z_3}{*} S_3) \underset{M}{*} (S_3 \underset{Z_2}{*} D_2)$, where $M$ is modular group whose finite presentation is $\langle D, C : D^2 = C^3 = 1 \rangle$. A general fragment of the coset diagram for the action of $\Gamma$ on $\mathbb{Q}(i)$ will look as follows.

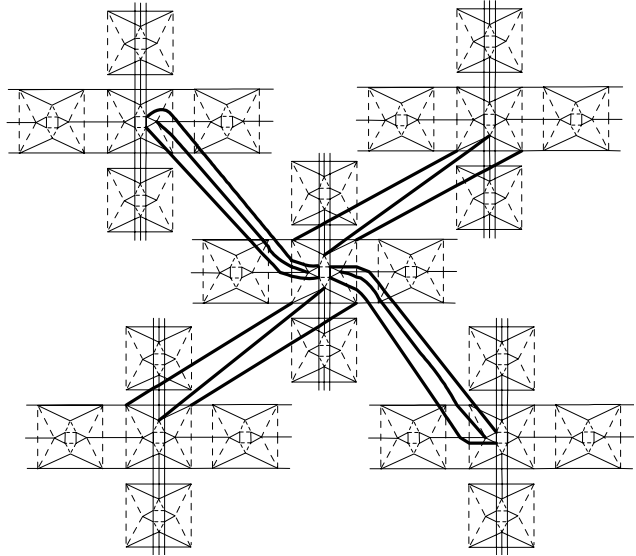

Fig.1

The generators $A$ and $C$ together form a Cayley's diagram of $A_4$. Two diagrams of $A_4$ are joined by edges of generator $D$. The edges have joined the 3 vertices of a triangle with broken edges of one diagram of $A_4$ to that of another diagram of $A_4$. The generators $A$ and $D$ together generate $S_3$. A triangle having broken edges is common between the diagram of $A_4$ and $S_3$. Algebraically, $\mathbb{Z}_3$ is an amalgam of $A_4$ and $S_3$. One diagram of $A_4$ can be joined to four other diagrams of $A_4$ by edges representing the generator $D$. We have called the fragment of the coset diagram which consists of diagrams of $A_4$ and edges of $D$ as a fragment $A - C - D$ because this fragment of coset diagram is composed by the use of generators $A, C,$ and $D$.

The two fragments $A - C - D$ have joined by bold edges representing the generator $B$ as shown in Fig. 1. The bold edges have joined the 3 vertices of triangle with unbroken edges of one diagram of $A_4$ of one fragment $A - C - D$ to that of another fragment $A - C - D$. The generators $B$ and $C$ together generate $S_3$. The generators $B$ and $D$ generate the group $D_2$. A bold edge is common between the diagram of $S_3$ and $D_2$. Algebraically, $\mathbb{Z}_2$ is an amalgam of $S_3$ and $D_2$.

Coset diagrams for the orbit of $\Gamma$ on biquadratic field $\mathbb{Q}(i,\sqrt{3})$ give some interesting information. We have shown that there is a finite number of ambiguous numbers in the orbit $\Gamma\alpha$, where $\alpha$ is ambiguous, and that they form a closed path and it is the only closed path in the orbit $\Gamma\alpha$. We have classified all the ambiguous numbers in the orbit. We need Propositions 10 to 9 to obtain one ambiguous number from the other.

**Proposition 1.** *If $\alpha = \frac{a+b\sqrt{3}}{c} \in \mathbb{Q}(i,\sqrt{3})$ is a totally positive real quadratic irrational number, then $C(\alpha)$ and $C^2(\alpha)$ are totally negative.*

**Proof.** Let $\alpha$ be a totally positive quadratic number. Then there are two possibilities either $a,d,c > 0$ or $a,d,c < 0$, where $d = \frac{a^2-3b^2}{c}$. When $a,d,c > 0$, then $C(\alpha) = C(\frac{a+b\sqrt{3}}{c}) = \frac{-a-d+b\sqrt{3}}{d}$.

Here $a_1 = -a-d < 0$, $c_1 = d > 0$, and $d_1 = 2a+c+d > 0$. This shows that $C(\alpha)$ is totally negative. Also

$$C^2(\alpha) = \frac{-1}{1+\alpha} = \frac{-a-c+b\sqrt{3}}{2a+c+d}.$$

Here $a_2 = -a-c < 0$, $c_2 = 2a+c+d > 0$, and $d_2 = \frac{a_2^2-3b^2}{c_2} = \frac{a^2+c^2+2ac-3b^2}{2a+c+d} = \frac{c(2a+c+d)}{2a+c+d} = c > 0$. This shows that $C^2(\alpha)$ is totally negative. Similarly, it can be proved that when $a,d,c < 0$, then $C(\alpha)$ and $C^2(\alpha)$ are totally negative. ∎

If a number $\alpha$ is ambiguous, then $\alpha\bar{\alpha} = \frac{a^2-3b^2}{c^2} < 0$, that is, $a^2-3b^2 < 0$. Thus, in other words $\alpha$ is ambiguous when $dc < 0$, where $d = \frac{a^2-3b^2}{c}$.

**Lemma 1.** *Transformations $B$ and $D$ map an ambiguous number to an ambiguous number.*

**Proof.** Let $\alpha = \frac{a+b\sqrt{3}}{c} \in \mathbb{Q}(i,\sqrt{3})$ be an ambiguous number, where $a,b,c \in \mathbb{Z}$. This implies that $\alpha\bar{\alpha} < 0$, that is, $\frac{a^2-3b^2}{c^2} < 0$ which further implies that $a^2-3b^2 < 0$. Now

$$B(\alpha) = \frac{1}{\alpha} = \frac{a-b\sqrt{3}}{d}, \text{ and}$$

$$D(\alpha) = \frac{-1}{\alpha} = \frac{a-b\sqrt{3}}{-d},$$

imply that $(B(\alpha))(\overline{B(\alpha)}) = \frac{a^2-3b^2}{d^2} = D(\alpha)(\overline{D(\alpha)})$.

Here $\frac{a^2-3b^2}{d^2} < 0$, since $a^2-3b^2 < 0$ and $d^2 > 0$. This shows that $B(\alpha)$ and $D(\alpha)$ are ambiguous numbers. ∎

**Remark 1.** *The value of $b$ is invariant for the elements of the form $\alpha = \frac{a+b\sqrt{3}}{c}$ in $\Gamma\alpha$, where $\alpha \in \mathbb{Q}(i,\sqrt{3})$.*

**Lemma 2.** *If $\alpha$ is not ambiguous, then so are $D(\alpha)$ and $B(\alpha)$.*

**Proof.** Suppose $\alpha = \frac{a+b\sqrt{3}}{c} \in \mathbb{Q}(i,\sqrt{3})$ is not ambiguous number, where $a,b,c \in \mathbb{Z}$. This means that $\alpha\bar{\alpha} \geq 0$, that is, $\frac{a^2-3b^2}{c^2} \geq 0$ which further implies that $a^2 - 3b^2 \geq 0$ because $c^2 \neq 0$. Since $B(\alpha) = \frac{1}{\alpha} = \frac{a-b\sqrt{3}}{d}$, and $D(\alpha) = \frac{-1}{\alpha} = \frac{-a+b\sqrt{3}}{d}$, so $(B(\alpha))(\overline{B(\alpha)}) = \frac{a^2-3b^2}{d^2} = D(\alpha)(\overline{D(\alpha)})$. Here $\frac{a^2-3b^2}{d^2} \geq 0$, because $a^2 - 3b^2 \geq 0$ and $d^2 > 0$. This shows that $B(\alpha)$ and $D(\alpha)$ are not ambiguous numbers. ∎

**Proposition 2.** *If $\alpha$ is an ambiguous number, then $A(\alpha)$ and $A^2(\alpha)$ are not ambiguous.*

**Proof.** Let $\alpha = \frac{a+b\sqrt{3}}{c}$ be an ambiguous number. This means that $\alpha\bar{\alpha} < 0$, that is, $\frac{a^2-3b^2}{c^2} < 0$. This implies that $a^2 - 3b^2 < 0$. After rationalization of $A(\alpha) = \frac{1}{\alpha-i} = \frac{c}{a+b\sqrt{3}-ci}$, the imaginary part is $(a^2c^2 + c^4 + 3b^2c^2 - 2abc^2\sqrt{3})$. This means that $A(\alpha)$ will be ambiguous if and only if $(a^2c^2 + c^4 + 3b^2c^2 - 2abc^2\sqrt{3})i = 0$. Since $i \neq 0$, it implies that $a^2c^2 + c^4 + 3b^2c^2 - 2abc^2\sqrt{3} = 0$. The term on left hand side will be zero only if $c = 0$. But $c$ cannot be zero because otherwise $\alpha = \infty$ will be not ambiguous. This shows that $A(\alpha)$ is not ambiguous.

Also, $A^2(\alpha) = \frac{1+i\alpha}{\alpha} = \frac{a-b\sqrt{3}+di}{d}$ will be ambiguous when imaginary part becomes zero, that is, $d = 0$. But $d$ cannot be zero because otherwise $A^2(\alpha)$ will become $\infty$. This shows that $A^2(\alpha)$ is not ambiguous. ∎

**Proposition 3.** *If $\alpha = \frac{a+b\sqrt{3}}{c} \in \mathbb{Q}(i,\sqrt{3})$ is an ambiguous number, then one of $C(\alpha)$ and $C^2(\alpha)$ is ambiguous and the other is totally negative.*

**Proof.** Suppose that $\alpha$ is a positive ambiguous number. Then by Proposition 1, the information can be tabulated as follows.

| $\alpha$ | $C(\alpha)$ | $C^2(\alpha)$ | $\bar{\alpha}$ | $\overline{C(\alpha)}$ | $\overline{C^2(\alpha)}$ |
|---|---|---|---|---|---|
| + | − | − | − | + | − |
|   |   |   | − | − | + |

Similarly, if $\alpha$ is a negative ambiguous number, then the information about $C(\alpha)$, $C^2(\alpha), \bar{\alpha}, \overline{C(\alpha)}$ and $\overline{C^2(\alpha)}$ can be tabulated as follows.

| $\alpha$ | $C(\alpha)$ | $C^2(\alpha)$ | $\bar{\alpha}$ | $\overline{C(\alpha)}$ | $\overline{C^2(\alpha)}$ |
|---|---|---|---|---|---|
| − | + | − | + | − | − |
| − | − | + |   |   |   |

Therefore, from the above tables it can easily be deduced that one of $C(\alpha)$ and $C^2(\alpha)$ is ambiguous and the other is totally negative. ∎

**Proposition 4.** *If* $\alpha = \frac{a+b\sqrt{3}}{c} \in \mathbb{Q}(i,\sqrt{3})$ *such that* $d = \frac{a^2-3b^2}{c}$ *is an integer, then the following hold:*

(i) $d$ *of* $C(\alpha)$ *and* $C^2(\alpha)$ *are integers,*

(ii) $d$ *of* $B(\alpha)$ *is an integer,*

(iii) $d$ *of* $D(\alpha)$ *is an integer.*

**Proof.** (i) Let $\alpha = \frac{a+b\sqrt{3}}{c} \in \mathbb{Q}(i,\sqrt{3})$ such that $d = \frac{a^2-3b^2}{c} \in \mathbb{Z}$. Then $C(\alpha) = \frac{-a-d+b\sqrt{3}}{d}$, where $d_1 = 2a+c+d \in \mathbb{Z}$, and $C^2(\alpha) = \frac{-a-c+b\sqrt{3}}{2a+c+d}$, where $d_2 = c \in \mathbb{Z}$.

(ii) Since $B(\alpha) = \frac{a-b\sqrt{3}}{d}$, the value of $d$ of $B(\alpha)$ is $\frac{a^2-3b^2}{d} = c \in \mathbb{Z}$.

(iii) Also $D(\alpha) = \frac{a-b\sqrt{3}}{-d}$ implies that $d$ of $D(\alpha)$ is $\frac{a^2-3b^2}{-d} = -c \in \mathbb{Z}$. ∎

A closed path whose all vertices are ambiguous numbers is called the closed path of ambiguous numbers. Lemma 3, Theorem 1 and Theorem 2 show that there is finite number of ambiguous numbers in one orbit and they form a single closed path in the orbit. In this single closed path the value of $b$ in $\frac{a+b\sqrt{3}}{c}$ remain invariant.

**Lemma 3.** *Ambiguous numbers in* $\Gamma\alpha$ *are finite.*

**Proof.** Let $\alpha = \frac{a+b\sqrt{3}}{c}$. It will be ambiguous when $\frac{a^2-3b^2}{c^2} < 0$ or $a^2 - 3b^2 < 0$ or $a^2 < 3b^2$. This shows that the values of $a$ is finite which satisfy the condition $a^2 < 3b^2$ for constant value of $b$. By Remark 1, the value of $b$ remain invariant for the numbers of the form $\frac{a+b\sqrt{3}}{c}$ in $\Gamma\alpha$, where $a,b,c \in \mathbb{Z}$. By Proposition 4, $d$ is integer, this implies that $c$ divides $(a^2 - 3b^2)$, so values of $c$ are also finite. As values of $a$ and $c$ are finite and value of $b$ is fixed in a orbit so ambiguous numbers of the form $\frac{a+b\sqrt{3}}{c}$ are also finite in an orbit. ∎

**Theorem 1.** *In a coset diagram for* $\Gamma\alpha$, *where* $\alpha$ *is an ambiguous number, the ambiguous numbers form a closed path.*

**Proof.** Suppose $\alpha_1$ is an ambiguous number in $\Gamma\alpha$. By $\alpha_1^k$ we mean ambiguous number $\alpha_1$ in $k^{\text{th}}$ triangle. Then by Proposition 3, either $C(\alpha_1^k)$ or $C^2(\alpha_1^k)$ is ambiguous. Let $C(\alpha_1^k) = \alpha_2^k$ be ambiguous. This means each triangle with unbroken edges in the coset diagram for $\Gamma$ contains two ambiguous numbers. By Proposition 10, $B(\alpha_1^k) = \alpha_1^m$, $B(\alpha_2^k) = \alpha_2^m$ and $D(\alpha_2^k) = \alpha_1^{k+1}, D(\alpha_1^k) = \alpha_2^{k-1}$ are also ambiguous. So within the $k^{\text{th}}$ triangle, the generators $D$ or $B$ are used to reach to the next ambiguous numbers as shown in Fig.2. Since $(BC)^2 = 1$, therefore there are just two paths, one from $k^{\text{th}}$ triangle and another from $m^{\text{th}}$ triangle, and generator $B$ connects these two paths.

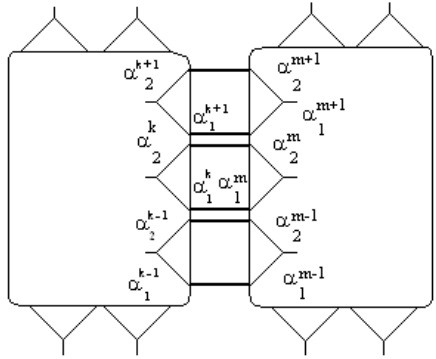

Fig.2

We can continue in this way as shown in Fig.2 and since by Lemma 3, there are only a finite number of ambiguous numbers, therefore after a finite number of steps we reach the vertex, $\alpha_1^{k+n} = \alpha_1^{k-1}$. Hence ambiguous numbers form a closed path in the coset diagram. ∎

**Theorem 2.** *If $\alpha$ is an ambiguous number, then in $\Gamma\alpha$ there is only one closed path of ambiguous numbers.*

**Proof.** In a diagram of $A_4$ there are four triangles with unbroken edges. So the diagram of $A_4$ in which $\alpha = \alpha_0$ is a vertex of one of the triangles with unbroken edges contains three more triangles with unbroken edges. According to Proposition 3, if $\alpha_0$ is an ambiguous number, then one of $C(\alpha_0)$ and $C^2(\alpha_0)$ is also an ambiguous number. Let us denote this ambiguous number by $\alpha_1$. Then by Lemma 1, these two ambiguous numbers are joined with other ambiguous numbers by generators $B$ and $D$. Therefore suppose $B(\alpha_0) = \alpha_2$, $B(\alpha_1) = \alpha_3$ and $D(\alpha_0) = \alpha_4$, $D(\alpha_1) = \alpha_5$.

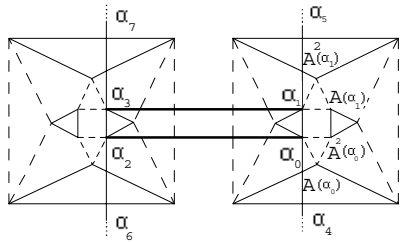

Fig.3

As the vertices of other three triangles with unbroken edges are also vertices of triangles with broken edges, according to Proposition 2 unambiguous, $A(\alpha_j)$ and $A^2(\alpha_j)$ are not ambiguous numbers but they contain $i$, where $j = 0, 1$. By applying transformations $C$ and $C^2$ on them, one gets not ambiguous numbers having $i$. So this diagram for $A_4$ contains only one triangle having two ambiguous numbers. If we expand this diagram and apply generators $B$ and $D$ on numbers which are not ambiguous, then by Lemma 2, no further ambiguous numbers are found. Since these numbers which are not ambiguous contain $i$, so by applying transformations $C$ and $C^2$ on them, one gets numbers which are not ambiguous but contain $i$. According to

Lemma 1 and Proposition 3, the generators $B, C$ and $D$ map $\alpha_0$ and $\alpha_1$ to other ambiguous numbers, namely, $\alpha_2, \alpha_3, \alpha_4, \ldots, \alpha_n$, since they are finite according to Lemma 3. By Theorem 1 they form a closed path. As one cannot found any further ambiguous numbers from the numbers which are not ambiguous, this means that there is only one closed path of ambiguous numbers in $\Gamma \alpha_0$. ∎

**Theorem 3.** *If $\alpha$ is an ambiguous number of the form $\frac{a+k\sqrt{3}}{c} \in \mathbb{Q}(i, \sqrt{3})$, where $a, c \in \mathbb{Z}$ and $k$ is a constant integer, then*
(i) $a^2 < 3k^2$,
(ii) $c$ *is a divisor of* $a^2 - 3k^2$.

**Proof.** ( i ) Let $\alpha = \frac{a+k\sqrt{3}}{c}$ be an ambiguous number. This means that $\alpha \bar{\alpha} < 0$ that is $\frac{a^2 - 3k^2}{c^2} < 0$. Since $c^2 > 0$, this implies that $a^2 - 3k^2 < 0$, which further implies that $a^2 < 3k^2$.

( ii ) The ambiguous number in the closed path containing $k\sqrt{3}$ are of the form $\frac{a+k\sqrt{3}}{c}$. We know that $d$ of $k\sqrt{3}$ is $-3k^2$. By Proposition 4, $d$ of other ambiguous numbers in the closed path containing $k\sqrt{3}$ are also integer. We have $d = \frac{a^2 - 3k^2}{c}$, it will be integer if and only if $c$ divides $a^2 - 3k^2$. This shows that $c$ is the divisor of $a^2 - 3k^2$. ∎

We denote the greatest common divisor of $a, b$ and $c$ by $(a, b, c)$, where $a, b, c \in \mathbb{Z}$. If $a$ divides $b$, then we denote it by $a \mid b$.

**Proposition 5.** *Let $m \neq k$ is a factor of $k \in \mathbb{Z}^+$ and $(a, k, c) \neq 1$. If $\alpha = \frac{á+m\sqrt{3}}{ć}$ such that $(á, m, ć) = 1$ and $ć \mid á^2 - 3m^2$, then ambiguous numbers of the form $\frac{a+m\sqrt{3}}{c}$ do not exist in the closed path of ambiguous numbers of the form $\frac{a+k\sqrt{3}}{c}$.*

**Proof.** Let $m \neq k$ be a factor of $k$. Let $\alpha = \frac{a+k\sqrt{3}}{c}$ be an ambiguous number such that $(a, k, c) \neq 1$. This means $\alpha$ can be written as $\frac{á+m\sqrt{3}}{ć}$ such that $(á, m, ć) = 1$. For $q \in \mathbb{Z}$ we have $áq = a, ćq = c$ and $mq = k$. Since $\alpha$ is ambiguous number, this means $á^2 < 3m^2$ and $ć \mid á^2 - 3m^2$. This means $\alpha$ is an ambiguous number of the form $\frac{á+m\sqrt{3}}{ć}$ and by Theorem 1, it forms a closed path of ambiguous numbers. By Remark 1, they are of the form $\frac{á+m\sqrt{3}}{ć}$. By Theorem one, it is the only closed path of ambiguous numbers of the form $\frac{á+m\sqrt{3}}{ć}$ in an orbit. So if $m \neq k$, then in the closed path all ambiguous numbers will be of the form $\frac{a+m\sqrt{3}}{c}$. Hence the ambiguous numbers of forms $\frac{a+k\sqrt{3}}{c}$ and $\frac{a+m\sqrt{3}}{c}$ lie in different closed paths. ∎

By Proposition 5, one can check whether an element of the form $\frac{a+k\sqrt{3}}{c}$, where $(a, k, c) \neq 1$, belongs to the closed path of ambiguous numbers of form $\frac{a+k\sqrt{3}}{c}$ or of the form $\frac{a+m\sqrt{3}}{c}$, where

$m \neq k$ is a factor of $k$. Let $\alpha = \frac{a+k\sqrt{3}}{c}$ such that $(a,k,c) \neq 1$, then $\alpha$ can be written as $\frac{á+m\sqrt{3}}{ć}$ such that $(á,m,ć) = 1$, where $m$ is a factor of $k$. If $ć \mid á^2 - 3m^2$, then $\alpha$ is of the form $\frac{a+m\sqrt{3}}{c}$ otherwise $\frac{a+k\sqrt{3}}{c}$. In the later case $c$ should divide $a^2 - 3k^2$. In other words if $\alpha = \frac{a+k\sqrt{3}}{c}$ is in its simplest form and $c \mid a^2 - 3k^2$, then $\alpha$ occurs in the closed path of the form $\frac{a+k\sqrt{3}}{c}$.

For instance let $k = 2$. The prime decomposition of $k$ is obviously $2 = 2 \times 1$. Therefore ambiguous numbers of the form $\frac{a+\sqrt{3}}{c}$ already exist in a closed path containing $\sqrt{3}$. So they will not exist in the closed path containing $2\sqrt{3}$. By using Theorem 3 and Proposition 5, it can be shown easily that the closed path containing $2\sqrt{3}$ has $32$ ambiguous numbers. For $\alpha = \frac{a+2\sqrt{3}}{c}$ to be ambiguous, $a^2 < 12$ which implies that $a = 0, \pm 1, \pm 2, \pm 3$. Let $d = \frac{a^2-12}{c}$. When $a = 0$, then $d = \frac{-12}{c}$ will be integer if and only if $c = \pm 1, \pm 2, \pm 3, \pm 4, \pm 6, \pm 12$. When $a = 0$ and $c = 2$, then $\alpha = \frac{2\sqrt{3}}{2} = \sqrt{3}$. Since $\sqrt{3}$ already exist in the closed path of ambiguous numbers of the form $\frac{a+\sqrt{3}}{c}$, we discard $c = \pm 2$ and its cofactor $\pm 6$. Also when $c = 6$, then $\alpha = \frac{2\sqrt{3}}{6} = \frac{\sqrt{3}}{3}$ and $3$ divides $3$. This shows that $\alpha$ belongs to the closed path of the form $\frac{a+\sqrt{3}}{c}$, therefore we also discard the cofactor $6$. Eventually, we have eight values of $c$, this means we have $8$ values of $\alpha$.

When $a = \pm 1$, then $d = \frac{1-12}{c} = \frac{-11}{c}$, $d \in \mathbb{Z}$ if and only if $c = \pm 1, \pm 11$. So we have four values of $c$. When $a = 2$, $d = \frac{4-12}{c} = \frac{-8}{c}$, $d \in \mathbb{Z}$ if and only if $c = \pm 1, \pm 2, \pm 4, \pm 8$. When $a = \pm 2$ and $c = \pm 2$, $\alpha = \frac{\pm 2 \pm 2\sqrt{3}}{2} = \pm 1 \pm \sqrt{3}$ is of the form $\frac{a+\sqrt{3}}{c}$. So we discard $2$ and $4$ (its cofactor). When $a = \pm 3$, then $d = \frac{9-12}{c} = \frac{-3}{c}$ implies that $c = \pm 1, \pm 3$. Here we have twelve values of $c$ for $a = \pm 1, \pm 2, \pm 3$. This means that there are $24$ values of $\alpha$ because of the positive and negative values of $a$. Thus the total ambiguous numbers of the form $\frac{a+2\sqrt{3}}{c}$ in an orbit of the coset diagram obtained by action of $\Gamma$ on $\mathbb{Q}(i, \sqrt{3})$ are $24 + 8 = 32$.

A fragment of the coset diagram obtained by action of $\Gamma$ on $\mathbb{Q}(i, \sqrt{3})$ containing closed path of ambiguous numbers of form $\frac{a+2\sqrt{3}}{c}$, is shown below. In Fig.4, here are two layers of fragment $A - C - D$. These two layers are connected by bold edges. The distinction between two layers are the conjugates, that is, if $\alpha = \frac{a+2\sqrt{3}}{c}$ occurs in one layer, then $\bar{\alpha} = \frac{a-2\sqrt{3}}{c}$ occurs in the second layer. Thus there is a mapping, say $S$, from one layer to the other, defined by $S: \alpha \mapsto \bar{\alpha}$. It can be noted here that in each diagram of $A_4$, there are exactly two ambiguous numbers.

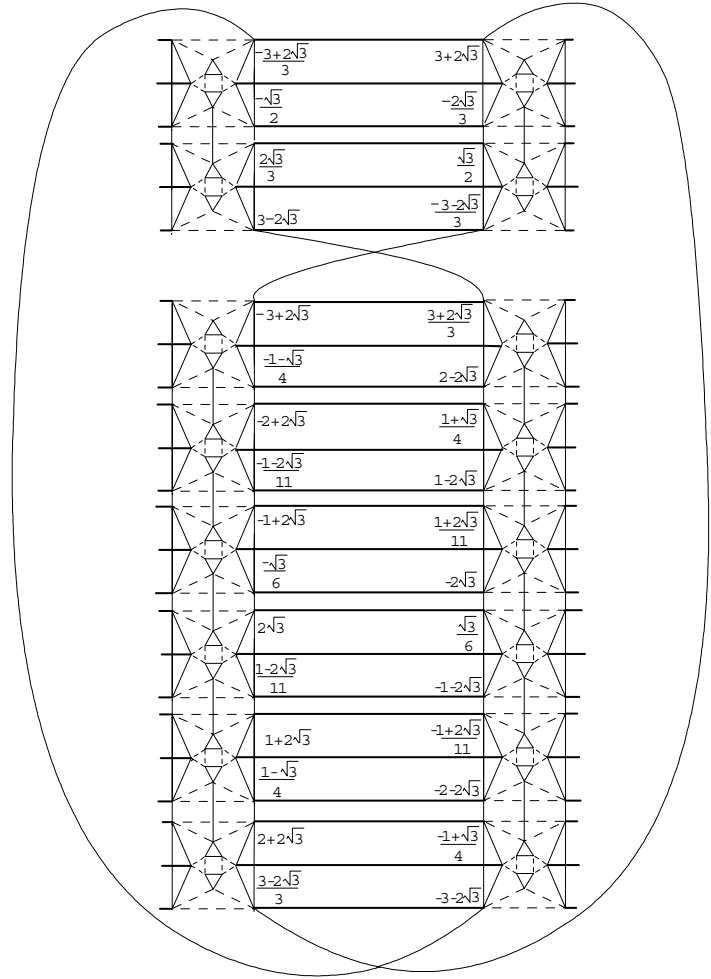

Fig.4

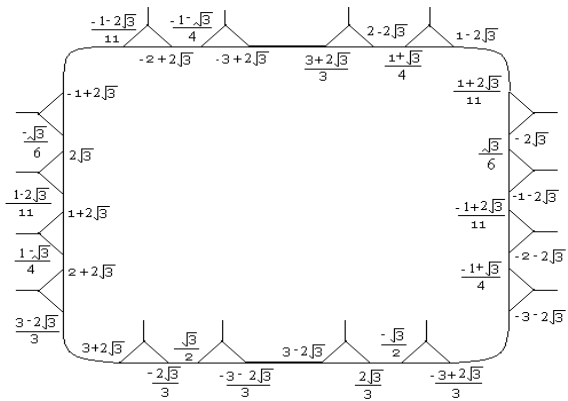

Fig.5

In this diagram we have shown only the closed path of ambiguous numbers of the form $\frac{a+2\sqrt{3}}{c}$ out of the fragment of Fig.4 . In Fig.5 , we apply repeatedly the generators $C, C^2$ and $D$ to obtain the ambiguous numbers in the path from $\frac{-3-2\sqrt{3}}{3}$ to $-3+2\sqrt{3}$ in one layer. The ambiguous numbers in the path from $\frac{3+2\sqrt{3}}{3}$ to $3-2\sqrt{3}$ are in second layer. These two layers are connected by bold edges. For instance, $-3+2\sqrt{3}$ and $\frac{3+2\sqrt{3}}{3}$ are connected by a bold edge as shown below. This is the only closed path of ambiguous numbers in the orbit $\Gamma(2\sqrt{3})$.